\documentclass[oneside]{amsart}
\usepackage{amsmath,amssymb,color,amsthm,graphicx,cite}
\usepackage[T1]{fontenc}
\usepackage[utf8]{inputenc}
\usepackage{color,bm}
\usepackage{enumerate} 
\usepackage[backref]{hyperref}
\RequirePackage{luatex85}
\usepackage[all]{xy}
\newtheorem{theorem}{Theorem}[section]
\newtheorem{proposition}[theorem]{Proposition}
\newtheorem{lemma}[theorem]{Lemma}

\theoremstyle{definition}
\newtheorem{definition}{Definition}
\newtheorem{main}{Theorem}

\newtheorem{main_cor}[main]{Corollary}

\def\F{\mathcal{F} }

\def\R{\mathbb{R} } 
 
\def\Z{\mathbb{Z} }

\def\R{\mathbb{R} }

\def\-{\ominus} 
\def\+{\oplus} 
\def\0{\circ}

\makeatletter
\@namedef{subjclassname@2020}{%
  \textup{2020} Mathematics Subject Classification}
\makeatother

\usepackage{color,bm}

\title[]{Filtrations Indexed by Attracting Levels and their Applications}


\author{Yusuke Imoto}
\address{Institute for the Advanced Study of Human Biology, Kyoto University, Yoshida Konoe-cho, Sakyo-ku, Kyoto, 606-8501, Japan}
\email{imoto.yusuke.4e@kyoto-u.ac.jp}

\author{Tomoo Yokoyama}
\address{Department of Mathematics, Faculty of Science, Saitama University, Shimo-Okubo 255, Sakura-ku, Saitama-shi, 338-8570 Japan}
\address{Department of Mathematics \& Statistics, McMaster University, 1280 Main Street West, Hamilton, Ontario, L8S 4L8, Canada\\}
\address{The Fields Institute for Research in Mathematical Sciences, 222 College Street, Toronto, Ontario, M5T 3J1, Canada}
\email{tyokoyama@rimath.saitama-u.ac.jp}

\begin{document}

\begin{abstract}
We introduce a new class of filtrations indexed by attracting levels in dynamical systems, providing novel inputs for persistent homology and related methods in topological data analysis. 
These filtrations quantify, in a forward direction, the sensitivity of trajectories with respect to attractors under perturbations and, in a backward direction, the perturbation magnitude at which attraction breaks down. 
The construction applies not only to maps on metric spaces but also to general partial maps with cost functions, yielding a filtration-theoretic framework with connections to algebraic topology. This generality ensures complementary filtrations when terminal states are good or bad, inducing natural decompositions of the underlying space. As an illustration, we apply the framework to ensemble forecasts of tropical cyclones, where the filtrations identify regions of heightened sensitivity, demonstrating the potential of our approach as a new tool for topological data analysis applied to dynamical systems.
\end{abstract}

\maketitle

\section{Introduction}

From the perspective of dynamical systems, Birkhoff introduced the concept of recurrent points and investigated the asymptotic behavior of orbits \cite{Birkhoff}.
The concept of $\varepsilon$-pseudo-orbits was introduced by Conley~\cite{conley1978isolated,conley1988gradient} and Bowen~\cite{bowen1975omega,Bowen1975,bowen1978axiom} for continuous and discrete dynamical systems, and has become a fundamental tool in the study of dynamical systems, particularly through numerical simulations.
Building on the concept of $\varepsilon$-pseudo-orbits, Conley also defined a weaker form of recurrence, known as chain recurrence \cite{conley1988gradient}.
In addition, the Conley theory says that dynamical systems on compact metric spaces can be decomposed into blocks,  each of which is either chain recurrent or gradient (cf. Conley\cite{conley1978isolated}).
%
In our previous study\cite{yokoyama2025coarse}, the concept of $\varepsilon$-recurrence was later introduced to analyze recurrent behaviors and persistence of attractors, and implies filtrations associated with dynamical systems that characterize the strength of recurrence properties.

Unifying these studies, this paper demonstrates that the idea behind this filtration can be applied to the analysis of systems with gradient-like dynamics or attractors, by using a concept similar to $\varepsilon$-pseudo-orbit. 
In fact, we construct filtrations associated with subsets in the complement of the domain of a mapping, and extend the construction to partial maps on sets equipped with cost functions, thereby obtaining a categorical and filtration-theoretic framework.

\begin{main}\label{main:01}
Let $(Y,d)$ be a metric space, $Y' \subset Y$ disjoint non-empty subsets, and $F \colon Y' \to Y$ a mapping. 
For any subset $A \subset Y - Y'$ with $A \cap \mathop{\mathrm{Im}}F \neq \emptyset$,  the families $(A_{F,\varepsilon})_{\varepsilon \in \R}$ {\rm(}see Definitions~\ref{def:pos_att}, \ref{def:zero_att}, and \ref{def:neg_att} below{\rm)}, and $(A_{F,\varepsilon_{\Sigma}})_{\varepsilon \in \R}$ {\rm(}see Definitions~\ref{def:pos_sum_att}, \ref{def:zero_sum_att}, and \ref{def:neg_sum_att} below{\rm)} with respect to $d$ are filtrations on $Y$. 
\end{main}

A more general statement is described in Theorem~\ref{thm:filtration} for any partial map on a set equipped with a cost function. 
The previous result essentially implies the following filtration, while its rigorous justification is provided by the general statement Theorem~\ref{thm:filtration}.

\begin{main_cor}\label{main_cor:01}
Let $f \colon X \to X$ be a mapping on a metric space $(X,d)$,  $a<b \in \Z_{\geq 0}$ integers, and $Y := \{a, a+1, \ldots , b \} \times X$.   
Define a mapping $F \colon Y - (\{b\} \times X) \to Y$ by $F(t,x):= (t+1, f(x))$ and a function $c_d \colon Y \times Y \to [0,\infty]$ as follows: 
\[
c_d ((t,x),(t',x')) := \begin{cases}
d(x,x') & \text{if } t=t'\\
\infty & \text{if } t \neq t' 
\end{cases}
\]
For any subset $A_X \subseteq X$ with $A_X \cap \mathop{\mathrm{Im}}f \neq \emptyset$, 
 the families $(A_{F,\varepsilon})_{\varepsilon \in \R}$ and $(A_{F,\varepsilon_{\Sigma}})_{\varepsilon \in \R}$ with respect to $c_d$ are filtrations, where $A := \{b\} \times A_X \subset Y$. 
\end{main_cor}

Notice that this filtration can measure the robustness of attractors in gradient-like dynamics and provide a lower bound on the magnitude of perturbations required for a orbit to no longer enter an attractor. 
Furthermore, it can give a lower bound on the energy needed to perturb the system so that it enters a different attractor, as well as a lower bound on the energy magnitude required at each step.

In addition, when the terminal state lies in either $G$ or $B$, the space decomposes into disjoint subsets: those that can be $\varepsilon$-controlled to reach $G$ and those that inevitably end up in $B$ (see Theorem~\ref{thm:eps_GBseparation} and Theorem~\ref{thm:eps_sum_GBseparation} for details).
This partition induces a natural filtration.
We also show that, for any $\varepsilon \in \R$, the corresponding subsets of dead-end states partition the space and induce a natural filtration (see Theorem~\ref{thm:eps_GBcovering} for details).

As an example of the proposed method, we introduce an application to meteorological data, specifically ensemble weather forecast data for a tropical cyclone \cite{oettli2024objective,oettli2025meteorological}, which is a crucial challenge in mitigating extreme weather events. 
{
By interpreting the parameter $\varepsilon$ as the magnitude of perturbations in initial conditions, we explore the times and locations at which small perturbations could have influenced the eventual trajectory of the tropical cyclone.
}

Moreover, the presence of a natural filtration in our framework suggests a direct connection with tools from topological data analysis (TDA), particularly persistent homology. 
While most filtrations studied in TDA arise from simplicial complexes or metric constructions, the filtrations introduced here originate from dynamical considerations and thus provide persistence theory with a new class of filtrations. 
This perspective indicates that persistence invariants can be used to detect structural transitions and robustness properties in dynamical systems. 
In particular, the proposed construction provides a new class of filtrations that open a natural avenue for applying TDA methods to problems of stability, sensitivity, and decomposition in dynamics, thereby bridging applied topology and dynamical systems.

The present paper consists of four sections.
In the next section, we recall and introduce some concepts of combinatorics and dynamical systems. 
Moreover, we describe fundamental properties and show Theorem~\ref{main:01} and a more general statement, Theorem~\ref{thm:filtration}. 
In Section~III, we construct filtrations for time extensions of partial maps. 
In addition, we demonstrated that when the final state falls into either a good or a bad state, the corresponding filtrations behave in complementary ways (see Theorems~\ref{thm:eps_GBseparation}--\ref{thm:eps_GBcovering} for details). 
In the final section, we discuss an application to meteorological data.

\section{Preliminaries}\label{sec:Preliminaries}

To analyze the various phenomena, we extend the framework of dynamical systems on metric spaces. 
As a preliminary step, we recall the definitions of filtration, partial map, and cost function as follows.

\begin{definition}\label{def:filtration}
Let $X$ be a set and $\F = \{ F_i \mid i \in I\} \subset 2^X$ a family indexed by a totally ordered set $I$, where $2^X$ is the power set of $X$. The family $\F$ is a {\bf filtration} if $X = \bigcup_{i \in I} F_i$ and $F_{i_1} \subseteq F_{i_2}$ for any pair $i_1 \leq i_2 \in I$.
\end{definition}

\begin{definition}\label{def:partial_map}
A map $f \colon X' \to Y$ from a subset $X' \subseteq X$ to a set $Y$ is called a {\bf partial map} from $X$ to $Y$ and denoted by $f \colon X \rightharpoonup Y$.
\end{definition}


\begin{definition}\label{def:positive_orbit}
For any partial map $f \colon X \rightharpoonup X$ and any point $x \in X$, the {\bf non-negative orbit} $O^{\geq 0}_f(x)$ is defined as follows:
\[
O^{\geq 0}_f(x) =
\begin{cases}
\{ x \} & \text{if } x \notin \mathop{\mathrm{dom}}f \\
\{ x, f(x), \ldots , f^n(x) \} & 
\text{if there is }n \in \Z_{>0} \text{ such that }f^n(x) \notin \mathop{\mathrm{dom}}f \\
& \text{and }  \{ x, f(x), \ldots , f^{n-1}(x) \} \subseteq  \mathop{\mathrm{dom}}f \\
\{ f^n(x) \mid n \in \Z_{\geq 0} \}  & \text{otherwise}
\end{cases}
\]
\end{definition}


\begin{definition}\label{def:cost}
Let $X$ be a set and $c \colon X^2 \to [0,\infty]$ a function.
The function $c$ is a {\bf cost function} if $c(x,x) = 0$ for any $x \in X$. 
\end{definition}

Note that we do not require transitivity in the previous definition to analyze various phenomena. 

\begin{definition}\label{def:metric-like}
A cost function $c \colon X^2 \to [0,\infty]$ is {\bf non-degenerate} if $c^{-1}(0) = \{ (x,x) \mid x \in X\}$. 
\end{definition}

Notice that a metric is a non-degenerate cost function. 
We generalize the conventional framework of dynamical systems on metric spaces, defined by iterations of maps, to a broader framework of dynamical systems based on iterations of partial maps on sets equipped with a cost function.


For sets $B,C$, the symbol $C - B$ is used instead of the set difference $C \setminus B$ when $B \subseteq C$.

From now on, in this section, let $F \colon Y \rightharpoonup Y$ be a partial map on a set $Y$ with a cost function $c \colon Y^2 \to [0,\infty]$.

\subsection{$\varepsilon$-controlled path}

We introduce the following concepts to analyze controls for partial maps on a set with a cost function. 

\begin{figure}[t]
\begin{center}
\includegraphics[scale=0.6]{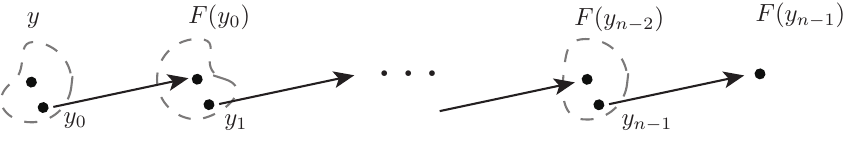}
\end{center} 
\caption{An $\varepsilon$-controlled path of length $n$ from the initial state $y$ to the terminal state $F(y_{n-1})$.}
\label{fig:chains}
\begin{center}
\includegraphics[scale=0.6]{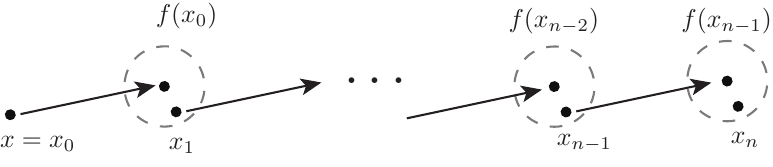}
\end{center} 
\caption{An $\varepsilon$-chain of length $n$ from a point $x$ to a point $x_n$.}
\label{fig:e-chains}
\end{figure}

\begin{definition}\label{def:e-p-chain}
Fix a number $\varepsilon \in [0,\infty]$, a point $y \in Y$, and a subset $A \subseteq Y$. 
Define the relation $y \overset{\exists}{\underset{F,\varepsilon, 0}{\rightharpoonup}} A$ as follows:
\[
y \overset{\exists}{\underset{F,\varepsilon, 0}{\rightharpoonup}} A \text{ if } y \in A
\]
Moreover, for any $n \in \Z_{> 0}$, the relation $y \overset{\exists}{\underset{F,\varepsilon, n}{\rightharpoonup}} A$ holds if there is a sequence $(y_0, \ldots , y_{n-1}) \in (\mathop{\mathrm{dom}}F)^{n}$ such that 
\[
c(y, y_0 ) \leq \varepsilon\]
\[
c(F(y_i), y_{i+1}) \leq \varepsilon\]
for any $i \in \{0,1, \ldots , n-2 \}$, and $F(y_{n-1}) \in A$, as shown in Figure~\ref{fig:chains}. 
Then $F(y_{n-1})$ is called a {{\rm(}\bf $\bm{\varepsilon}$-{\rm)}terminal state} of $y \overset{\exists}{\underset{F,\varepsilon,n}{\rightharpoonup}} A$, and the sequence 
\[
(y;y_0, \ldots , y_{n-1};F(y_{n-1}))
\]
is called an {\bf $\bm{\varepsilon}$-controlled path} of length $n$ from the initial state $y$ to the terminal state $F(y_{n-1})$.
Moreover, the relation $y \overset{\exists}{\underset{F,\varepsilon}{\rightharpoonup}} A$ holds if there is a non-negative number $n \in \Z_{\geq 0}$ such that $y \overset{\exists}{\underset{F,\varepsilon, n}{\rightharpoonup}} A$. 
\end{definition}

In the previous definition, the concept of $\varepsilon$-controlled path 
is similar to one of the $\varepsilon$-chains. 
However, it is defined differently from the following perspective of control:
To ensure that the current state does not become ``bad'' even instantaneously during control, the control is applied before the partial map, rather than after. This corresponds to applying control first, followed by time evolution. 
This is in contrast to the concept of an $\varepsilon$-chain as shown in Figure~\ref{fig:e-chains}, in which control is applied only after the partial map has been executed.
The time evolution does not end in the controlled state; rather, in correspondence with the fact that time continues to evolve afterward, it ends after the controlled state has been mapped.
We consider the set of points from which there is a $\varepsilon$-controlled path from a point as follows. 

\begin{definition}
For any $\varepsilon \in [0,\infty]$ and any point $y \in Y$, denote by $[y]_F^\varepsilon$ the set of points to which there is a $\varepsilon$-controlled path from $y$, called the $\varepsilon$-reachable set for $y$.
\end{definition}

In other words, we have that $[y]_F^\varepsilon = \{ y' \in Y \mid y \overset{\exists}{\underset{F,\varepsilon}{\rightharpoonup}} \{y' \} \}$. 
We have the following observation. 

\begin{lemma}\label{lem:pos_orbit}
The following statements hold for any $y \in Y$:
\\
{\rm(1)} If the cost function is non-degenerate, then $[y]_F^0 = O^{\geq 0}_F(y)$. 
\\
{\rm(2)} For any $\varepsilon_1 \leq \varepsilon_2 \in [0,\infty]$, we obtain $[y]_F^{\varepsilon_1} \subseteq [y]_F^{\varepsilon_2}$. 
\end{lemma}
%


\subsubsection{$\varepsilon$-attracting basin}

We generalize attracting basins as follows. 

\begin{definition}\label{def:pos_att}
Let $\varepsilon \in [0,\infty]$ be a number. 
For any subset $A \subseteq Y$, define the {\bf $\bm{\varepsilon}$-attracting basin} $A_{F,\varepsilon}$ of $A$ with respect to the cost $c$ as follows:
\[
A_{F,\varepsilon} := \{ y \in Y \mid y \overset{\exists}{\underset{F,\varepsilon}{\rightharpoonup}} A \} = \{ y \in Y \mid [y]_F^\varepsilon \cap A \neq \emptyset \} 
\]
\end{definition}

 
Notice that the $\varepsilon$-attracting basin is similar to the dual concept of the $(\varepsilon,t)$-approximation of the $\Omega$-limit set $\Omega'(C,\varepsilon,t,f)$ in \cite{conley1988gradient}.
As mentioned after Definition~\ref{def:e-p-chain}, the operations of partial map and perturbation are applied in the reverse order.
In addition, the approximation $\Omega'(C,\varepsilon,t,f)$ is a set obtained by perturbing points of a subset $C$, whereas, conversely, $A_{F,\varepsilon}$ is the set of points that can reach $A$ via a perturbation.

\subsubsection{Dead-end terminal states and $-\varepsilon$-attracting basin}

We introduce a concept corresponding to the ``depth within the interior of an attracting basin''.
To construct the concept, we introduce the following technical concept.

\begin{definition}
For any subset $A \subseteq Y$ and any point $y \in A_{F,0}$, the relation $y \overset{\forall}{\underset{F,0}{\rightharpoondown}} A$ holds if $\emptyset \neq [y]_F^0 \setminus \mathop{\mathrm{dom}}F \subseteq A$. 
\end{definition}

We obtain the following equivalence.

\begin{lemma}\label{lem:equv_zero}
If the cost function $c$ is non-degenerate, then the following statements are equivalent for any subset $A \subseteq Y$ and any point $y \in Y$:
\\
{\rm(1)} $y \overset{\forall}{\underset{F,0}{\rightharpoondown}} A$ {\rm(i.e.} $\emptyset \neq [y]_F^0 \setminus \mathop{\mathrm{dom}}F \subseteq A${\rm)}. 
\\
{\rm(2)} $y \overset{\exists}{\underset{F,0}{\rightharpoonup}} A \setminus \mathop{\mathrm{dom}}F$  {\rm(i.e.} $([y]_F^0 \setminus \mathop{\mathrm{dom}}F) \cap A \neq \emptyset${\rm)}. 
\end{lemma}

\begin{proof}
By definitions, assertion (1) implies assertion (2). 
Suppose that the cost function $c$ is non-degenerate and assertion (2) holds for the point $y \in Y$. 
By Lemma~\ref{lem:pos_orbit}, the non-degeneracy of $c$ implies that $[y]_F^0 = O^{\geq 0}_F(y)$. 
The following statements are equivalent for any $y \in Y$:
\[
\begin{split}
&([y]_F^0 \setminus \mathop{\mathrm{dom}}F) \cap A \neq \emptyset
\\
\Longleftrightarrow &  (O^{\geq 0}_F(y) \setminus \mathop{\mathrm{dom}}F) \cap A \neq \emptyset
\\
\Longleftrightarrow &  (O^{\geq 0}_F(y) \setminus \mathop{\mathrm{dom}}F) \cap (A  \setminus \mathop{\mathrm{dom}}F) \neq \emptyset
\\
\Longleftrightarrow &  \vert (O^{\geq 0}_F(y) \setminus \mathop{\mathrm{dom}}F) \cap (A  \setminus \mathop{\mathrm{dom}}F) \vert = 1
\\
\Longleftrightarrow &  \emptyset \neq O^{\geq 0}_F(y) \setminus \mathop{\mathrm{dom}}F \subseteq (A  \setminus \mathop{\mathrm{dom}}F) 
\\
\Longleftrightarrow &  \emptyset \neq [y]_F^0 \setminus \mathop{\mathrm{dom}}F \subseteq (A  \setminus \mathop{\mathrm{dom}}F) 
\end{split}
\]
This means that assertion (2) implies assertion (1). 
\end{proof}

We introduce the following concepts.

\begin{definition}\label{def:zero_att}
For any subset $A \subseteq Y$, define the {\bf $\bm{-0}$-attracting basin} $A_{F,-0}$ of $A$ with respect to the cost $c$ as follows:
\[
\begin{split}
A_{F,-0} := \{ y \in Y \mid y \overset{\forall}{\underset{F,0}{\rightharpoondown}} A \} &= \{ y \in A_{F,0} \mid y \overset{\forall}{\underset{F,0}{\rightharpoondown}} A \} 
\\
&= \{ y \in A_{F,0} \mid \emptyset \neq [y]_F^0 \setminus \mathop{\mathrm{dom}}F \subseteq A \} 
\end{split}
\]
\end{definition}

By definitions of $A_{F,0}$ and $A_{F,-0}$, we have $A_{F,-0} \subseteq A_{F,0}$ for any partial map $F \colon Y \rightharpoonup Y$ and any $A \subseteq Y$. 

\begin{definition}
Let $F \colon Y \rightharpoonup Y$ be a partial map and $\varepsilon \in (0,\infty]$ a number. 
For any subset $A \subseteq Y$ and any point $y \in A_{F,-0}$, the relation $y \overset{\forall}{\underset{F,\varepsilon}{\rightharpoondown}} A$ holds if $\emptyset \neq [y]_F^\varepsilon \setminus \mathop{\mathrm{dom}}F \subseteq A$. 
Then an element of $[y]_F^\varepsilon \setminus \mathop{\mathrm{dom}}F$ is called a {\bf {\rm(}dead-end{\rm)} ${\rm(}\bm{\varepsilon}$-{\rm)}terminal state} of $y \overset{\forall}{\underset{F,\varepsilon}{\rightharpoondown}} A$. 
\end{definition}

Notice that $[y]_F^\varepsilon \setminus \mathop{\mathrm{dom}}F \subseteq \{ y \} \cup  \mathop{\mathrm{Im}}F$ for any partial map  $F \colon Y \rightharpoonup Y$ and any $\varepsilon \in (0,\infty]$. 

\begin{definition}\label{def:neg_att}
Let $\varepsilon \in (0,\infty]$ be a number. 
For any subset $A \subseteq Y$, define the {\bf $\bm{-\varepsilon}$-attracting basin} $A_{F,-\varepsilon}$ of $A$ with respect to the cost $c$ as follows:
\[
\begin{split}
A_{F,-\varepsilon} := &\,  \{ y \in A_{F,-0} \mid y \overset{\forall}{\underset{F,\varepsilon}{\rightharpoondown}} A \}
\\
= &\,  \{ y \in A_{F,-0} \mid \emptyset \neq [y]_F^\varepsilon \setminus \mathop{\mathrm{dom}}F \subseteq A \}
\\
= &\,  \{ y \in Y \mid \emptyset \neq [y]_F^\varepsilon \setminus \mathop{\mathrm{dom}}F \subseteq A \}
\end{split}
\]
\end{definition}

Moreover, we have the following observation.

\begin{proposition}
The following statement holds for any $\varepsilon_1 \geq \varepsilon_2 \in [0,\infty]$:
\[
A_{F,-\varepsilon_1} \subseteq A_{F,-\varepsilon_2} \subseteq A_{F,-0} \subseteq A_{F,0} \subseteq A_{F,\varepsilon_2} \subseteq A_{F,\varepsilon_1}
\]
\end{proposition}

\subsection{$\sum$-perturbations}

In this subsection, we introduce the definitions corresponding to the case where ``the total energy is $\varepsilon$'', in contrast to the previous definitions where ``the maximum energy of $\varepsilon$ was assigned to each step''.

\begin{definition}
Fix a number $\varepsilon \in [0,\infty]$, a point $y \in Y$, and a subset $A \subseteq Y$. 
Define the relation $y \overset{\exists}{\underset{F,\varepsilon_{\Sigma}, 0}{\rightharpoonup}} A$ as follows:
\[
y \overset{\exists}{\underset{F,\varepsilon_{\Sigma}, 0}{\rightharpoonup}} A \text{ if } y \in A
\]
Moreover, for any $n \in \Z_{> 0}$, the relation $y \overset{\exists}{\underset{F,\varepsilon_{\Sigma}, n}{\rightharpoonup}} A$ holds if there are a sequence $(y_0, \ldots , y_{n-1}) \in (\mathop{\mathrm{dom}}F)^{n}$ and a sequence $(\varepsilon_0, \ldots , \varepsilon_{n-1}) \in [0,\infty]^{n}$ with $\varepsilon \geq \sum_{i=0}^{n-1} \varepsilon_i$ such that 
\[
c(y, y_0 ) \leq \varepsilon_0 \]
\[
c(F(y_i), y_{i+1}) \leq \varepsilon_{i+1}\]
for any $i \in \{0,1, \ldots , n-2 \}$, and $F(y_{n-1}) \in A$. 
Then $F(y_{n-1})$ is called a {\bf ${\rm(}\bm{\varepsilon_{\Sigma}}$-{\rm)}terminal state} of $y \overset{\exists}{\underset{F,\varepsilon_{\Sigma}, n}{\rightharpoonup}} A$, and the sequence $(y;y_0, \ldots , y_{n-1};F(y_{n-1}))$ is called an {\bf $\bm{\varepsilon_{\Sigma}}$-controlled path} of length $n$ from the initial state $y$ to the terminal state $F(y_{n-1})$.
Moreover, the relation $y \overset{\exists}{\underset{F,\varepsilon_{\Sigma}}{\rightharpoonup}} A$ holds if there is a non-negative number $n \in Z_{\geq 0}$ such that $y \overset{\exists}{\underset{F,\varepsilon_{\Sigma}, n}{\rightharpoonup}} A$. 
\end{definition}

\begin{definition}
For any $\varepsilon \in [0,\infty]$ and any point $y \in Y$, denote by $[y]_F^{\varepsilon_{\Sigma}}$ the set of points to which  there is a $\varepsilon_{\Sigma}$-controlled path from $y$.
\end{definition}

In other words, we have that $[y]_F^{\varepsilon_{\Sigma}} = \{ y' \in Y \mid y \overset{\exists}{\underset{F,\varepsilon_\Sigma}{\rightharpoonup}} \{y' \} \}$. 
We have the following observation. 

\begin{lemma}\label{lem:pos_sum_orbit}
The following statements hold for any $y \in Y$:
\\
{\rm(1)} If the cost function is non-degenerate, then $[y]_F^{0_\Sigma} = O^{\geq 0}_F(y) = [y]_F^{0}$. 
\\
{\rm(2)} For any $\varepsilon_1 \leq \varepsilon_2 \in [0,\infty]$, we obtain $[y]_F^{\varepsilon_{1 \Sigma}} \subseteq [y]_F^{\varepsilon_{2 \Sigma}}$. 
\end{lemma}


\subsubsection{$\varepsilon_{\Sigma}$-attracting basin}

We generalize attracting basins as follows. 

\begin{definition}\label{def:pos_sum_att}
Let $F \colon Y \rightharpoonup Y$ be a partial map and $\varepsilon \in [0,\infty]$. 
For any subset $A \subseteq Y$, define the {\bf $\bm{\varepsilon_{\Sigma}}$-attracting basin} $A_{F,\varepsilon_{\Sigma}}$ of $A$ with respect to the cost $c$ as follows:
\[
A_{F,\varepsilon_{\Sigma}} := \{ y \in Y \mid y \overset{\exists}{\underset{F,\varepsilon_{\Sigma}}{\rightharpoonup}} A \} = \{ y \in Y \mid [y]_F^{\varepsilon_{\Sigma}} \cap A \neq \emptyset \} 
\]
\end{definition}

By definition, we have the following observation.

\begin{lemma}\label{lem:two_inclusions}
For any $\varepsilon \in [0,\infty]$, we obtain $A_{F,\varepsilon_{\Sigma}} \subseteq A_{F,\varepsilon}$. 
\end{lemma}

\subsubsection{Dead-end $\varepsilon_{\Sigma}$-terminal states and $-\varepsilon_{\Sigma}$-attracting basin}

We also introduce a concept corresponding to the ``depth within the interior of an attracting basin''.
To construct the concept, we introduce the following technical concept.

\begin{definition}
For any subset $A \subseteq Y$ and any point $y \in A_{F,0_{\Sigma}}$, the relation $y \overset{\forall}{\underset{F,0_{\Sigma}}{\rightharpoondown}} A$ holds if $\emptyset \neq [y]_F^{0_\Sigma} \setminus \mathop{\mathrm{dom}}F \subseteq A$. 
\end{definition}

By the same argument of the proof of Lemma~\ref{lem:equv_zero}, replacing Lemma~\ref{lem:pos_orbit} with Lemma~\ref{lem:pos_sum_orbit} in the proof, we obtain the following observation.

\begin{lemma}\label{lem:equv_zero_sum}
If the cost function $c$ is non-degenerate, then the following statements are equivalent for any subset $A \subseteq Y$ and any point $y \in Y$:
\\
{\rm(1)} $y \overset{\forall}{\underset{F,0_\Sigma}{\rightharpoondown}} A$ {\rm(i.e.} $\emptyset \neq [y]_F^{0_\Sigma} \setminus \mathop{\mathrm{dom}}F \subseteq A${\rm)}. 
\\
{\rm(2)} $y \overset{\exists}{\underset{F,0_\Sigma}{\rightharpoonup}} A \setminus \mathop{\mathrm{dom}}F$ {\rm(i.e.} $([y]_F^{0_\Sigma} \setminus \mathop{\mathrm{dom}}F) \cap A \neq \emptyset${\rm)}. 
\end{lemma}

We introduce the following concepts.

\begin{definition}\label{def:zero_sum_att}
For any subset $A \subseteq Y$, define the {\bf $\bm{-0_{\Sigma}}$-attracting basin} $A_{F,-0_{\Sigma}}$ of $A$ with respect to the cost $c$ as follows:
\[
\begin{split}
A_{F,-0_\Sigma} := & \{ y \in Y \mid y \overset{\forall}{\underset{F,0_\Sigma}{\rightharpoondown}} A \} 
\\
= & \{ y \in A_{F,0_\Sigma} \mid y \overset{\forall}{\underset{F,0_\Sigma}{\rightharpoondown}} A \} 
\\
= & \{ y \in A_{F,0_\Sigma} \mid \emptyset \neq [y]_F^{0_\Sigma} \setminus \mathop{\mathrm{dom}}F \subseteq A \} 
\end{split}
\]
\end{definition}


\begin{definition}
Let $\varepsilon \in (0,\infty]$ be a number. 
For any subset $A \subseteq Y$ and any point $y \in A_{F,-0_{\Sigma}}$, the relation $y \overset{\forall}{\underset{F,\varepsilon_{\Sigma}}{\rightharpoondown}} A$ holds if $\emptyset \neq [y]_F^{\varepsilon_\Sigma} \setminus \mathop{\mathrm{dom}}F \subseteq A$. 
Then an element of $[y]_F^{\varepsilon_\Sigma} \setminus \mathop{\mathrm{dom}}F$ is called a {\bf {\rm(}dead-end{\rm)} {\rm(}$\bm{\varepsilon_{\Sigma}}$-{\rm)}terminal state} of $y \overset{\forall}{\underset{F,\varepsilon_{\Sigma}}{\rightharpoondown}} A$.
\end{definition}

\begin{definition}\label{def:neg_sum_att}
Let $\varepsilon \in (0,\infty]$ be a number. 
For any subset $A \subseteq Y$, 
 define the {\bf $\bm{-\varepsilon_{\Sigma}}$-attracting basin} $A_{F,-\varepsilon_{\Sigma}}$ of $A$ with respect to the cost $c$ as follows:
\[
\begin{split}
A_{F,-\varepsilon_\Sigma} := &\,  \{ y \in A_{F,-0_\Sigma} \mid y \overset{\forall}{\underset{F,\varepsilon_\Sigma}{\rightharpoondown}} A \}
\\
= &\,  \{ y \in A_{F,-0_\Sigma} \mid \emptyset \neq [y]_F^{\varepsilon_\Sigma} \setminus \mathop{\mathrm{dom}}F \subseteq A \}
\\
= &\,  \{ y \in Y \mid \emptyset \neq [y]_F^{\varepsilon_\Sigma} \setminus \mathop{\mathrm{dom}}F \subseteq A \}
\end{split}
\]
\end{definition}

We have the following observation.

\begin{lemma}\label{lem:023}
The following statements hold for any partial map $F \colon Y \rightharpoonup Y$, any $A \subseteq Y$, and any $\varepsilon \in [0,\infty]$:
\\
{\rm(1)} $A_{F,-\varepsilon} \subseteq A_{F,-\varepsilon_{\Sigma}} \subseteq A_{F,-0_{\Sigma}} = A_{F,-0} \subseteq A_{F,0} = A_{F,0_{\Sigma}} \subseteq  A_{F,\varepsilon_{\Sigma}} \subseteq A_{F,\varepsilon}$. 
\\
{\rm(2)} If the cost function $c$ is non-degenerate and $A \subseteq Y - \mathop{\mathrm{dom}}F$, then $A_{F,0_{\Sigma}} = A_{F,-0} = A_{F,0} = A_{F,-0_{\Sigma}}$. 
\end{lemma}

\begin{proof}
By definitions of $A_{F,\varepsilon_{\Sigma}}$ and $A_{F,-\varepsilon_{\Sigma}}$, assertion (1) holds. 
Suppose that the cost function $c$ is non-degenerate and $A \subseteq Y - \mathop{\mathrm{dom}}F$. 
By Lemma~\ref{lem:equv_zero} and Lemma~\ref{lem:equv_zero_sum}, assertion (2) follows from the above equivalence and assertion (1). 
\end{proof}

\subsection{Properties of increasing families}

We have the following observation.

\begin{lemma}\label{lem:inclusions}
For any $\varepsilon_1 \geq \varepsilon_2 \in [0,\infty]$, we have $A_{F,\varepsilon_{1\sum}} \supseteq A_{F,\varepsilon_{2\sum}}$ and $A_{F,-\varepsilon_{1\sum}} \subseteq A_{F,-\varepsilon_{2\sum}}$.
\end{lemma}



We have the following filtrations associated with any subset $A \subseteq Y$ with $A \cap \mathop{\mathrm{Im}} F \neq \emptyset$.

\begin{theorem}\label{thm:filtration}
Let $F \colon Y \rightharpoonup Y$ be a partial map on a set $Y$ with a cost function $c \colon Y^2 \to [0,\infty]$.
We have the following statements for any subset $A \subseteq Y$ with $A \cap \mathop{\mathrm{Im}} F \neq \emptyset$: 
\\
{\rm(1)} The families $(A_{F,\varepsilon})_{\varepsilon \in (-\infty,-0] \sqcup [0,\infty]}$ and $(A_{F,\varepsilon_{\Sigma}})_{\varepsilon \in (-\infty,-0] \sqcup [0,\infty]}$ are filtration. 
\\
{\rm(2)} If the cost function $c$ is non-degenerate and $A \subseteq Y - \mathop{\mathrm{dom}}F$, then the families $(A_{F,\varepsilon})_{\varepsilon \in \R \sqcup \{\infty\}}$ and $(A_{F,\varepsilon_{\Sigma}})_{\varepsilon \in \R \sqcup \{\infty\}}$ are filtrations. 
\\
{\rm(3)} If $c^{-1}(\infty) = \emptyset$, then $(A_{F,\varepsilon})_{\varepsilon \in (-\infty,-0] \sqcup [0,\infty)}$ and $(A_{F,\varepsilon_{\Sigma}})_{\varepsilon \in (-\infty,-0] \sqcup [0,\infty)}$ are filtrations. 
\\
{\rm(4)} If the cost function $c$ is non-degenerate and $A \subseteq Y - \mathop{\mathrm{dom}}F$ and $c^{-1}(\infty) = \emptyset$, then $(A_{F,\varepsilon})_{\varepsilon \in \R}$ and $(A_{F,\varepsilon_{\Sigma}})_{\varepsilon \in \R}$ are filtrations. 
\end{theorem}

\begin{proof}
From Lemma~\ref{lem:023}, assertion (1) (resp. (3)) implies assertion (2) (resp. (4)). 
By Lemma~\ref{lem:inclusions}, it suffices to show that $\bigcup_{\varepsilon \in \R_{\geq 0} \sqcup \{\infty\}} A_{F,\varepsilon} = Y$. 
From definitions of $(A_{F,\varepsilon})_{\varepsilon \in \R_{\geq 0}}$ and $(A_{F,\varepsilon_{\Sigma}})_{\varepsilon \in \R_{\geq 0}}$, if $(A_{F,\varepsilon_{\Sigma}})_{\varepsilon \in \R_{\geq 0}}$ is a filtration, then so is $(A_{F,\varepsilon})_{\varepsilon \in \R_{\geq 0}}$. 
Fix a point $y \in Y$. 
From $A \cap \mathop{\mathrm{Im}} F \neq \emptyset$, there is a point $a \in Y$ such that $F(a) \in A$. 
Then $y \overset{\exists}{\underset{F,c(y,a)_{\Sigma}, 1}{\rightharpoonup}} A$. 
This means that $y \in A_{F,c(y,a)_{\Sigma}}$. 
Therefore, we have $\bigcup_{\varepsilon \in \R_{\geq 0} \sqcup \{\infty\}} A_{F,\varepsilon_{\Sigma}} = Y$. 

Suppose that $c^{-1}(\infty) = \emptyset$. 
Then $\bigcup_{\varepsilon \in \R_{\geq 0}} A_{F,\varepsilon_{\Sigma}} = Y$. 
\end{proof}

Theorem~\ref{main:01} and Corollary~\ref{main_cor:01} follow from the previous theorem. 
Notice that the assumption $A \cap \mathop{\mathrm{Im}} F \neq \emptyset$ in the previous theorem is necessary. 
In fact, there is a map $F \colon \{0,1,2\} \to \{0,1,2\}$ defined by $0 = F(0) = F(1)$ and $1 = F(2)$, and put $A := \{2 \}$ such that $(A_{F,\varepsilon})_{\varepsilon \in \R}$ is not a filtration for any cost function on $\{0,1,2\}$. 
Indeed, since the image $\mathop{\mathrm{Im}} F = \{0,1 \}$ does not contain $2$, for any cost function on $\{0,1,2\}$ and for any $\varepsilon \in \R$, we have  $A_{F,\varepsilon} = \{ 2 \}$. 

\subsubsection{Debut functions}

We define the following functions. 

\begin{definition}\label{def:debut}
Define functions \(\underline{\varepsilon}, \underline{\varepsilon}_{\Sigma} \colon Y \times 2^Y \rightarrow \mathbb{R}\sqcup \{ \infty \}\) as follows:
\[
\underline{\varepsilon}(y; A) := \inf \{\varepsilon\in\R \mid y \in A_{F,\varepsilon}\},
\]
\[
\underline{\varepsilon}_{\Sigma}(y; A) := \inf \{\varepsilon\in\R \mid y \in A_{F,\varepsilon_{\Sigma}}\}
\] 
We call functions $\underline{\varepsilon}(\cdot ; A)$ and  $\underline{\varepsilon}_{\Sigma}(\cdot ; A)$ {\bf debut functions} with respect to $A$, because these function are analogous to the debuts in \cite{dellacherie1978probabilities}. 
\end{definition}

By the previous theorem, notice that the images of the restrictions of $\underline{\varepsilon}, \underline{\varepsilon}_{\Sigma}$ to the subset $Y \times  \mathop{\mathrm{Im}F}$ are contained in $\R$.
By definitions, we have the following observations.

\begin{lemma}\label{lem:debut_ineq1}
The following statements hold:
\\
{\rm(1)}
For any $y \not\in A_{F,0}$, we have \(0 \leq \underline{\varepsilon}(y; A) \leq \underline{\varepsilon}_{\Sigma}(y; A)
\).
\\
{\rm(2)} For any $y \in A_{F,0}$, we have \(0 \geq \underline{\varepsilon}(y; A) \geq \underline{\varepsilon}_{\Sigma}(y; A)\).
\end{lemma}

\begin{lemma}\label{lem:debut_ineq2}
The following statements hold for any $y \in Y$:
\\
{\rm(1)}
If $F(y) \not\in A_{F,0}$, then $\underline{\varepsilon}(y; A) \leq \underline{\varepsilon}(F(y); A)$ and $\underline{\varepsilon}_{\Sigma}(y; A) \leq \underline{\varepsilon}_{\Sigma}(F(y); A)$. 
\\
{\rm(2)} If $F(y) \in A_{F,0}$, then $0 \geq \underline{\varepsilon}(y; A) \geq \underline{\varepsilon}(F(y); A)$ and $0 \geq \underline{\varepsilon}_{\Sigma}(y; A) \geq \underline{\varepsilon}_{\Sigma}(F(y); A)$.
\end{lemma}

\subsection{Terminal states which need not be dead-end}

To characterize the $\pm\varepsilon$-($\Sigma$-)attracting basins of discrete dynamical systems that contain ``dead ends'', we introduce the following concepts, which will be used in the next section.

\begin{definition}
Let $\varepsilon \in (0,\infty]$ be a number. 
For any point $y \in Y$ and any subset $A \subseteq Y$ and for any $n \in \Z_{\geq 0}$, the relation $y \overset{\forall}{\underset{F,\varepsilon, n}{\rightharpoonup}} A$ holds if $\emptyset \neq \{ y' \mid y \overset{\exists}{\underset{F,\varepsilon, n}{\rightharpoonup}} \{y' \} \} \subseteq A$. 
Then any element of $\{ y' \mid y \overset{\exists}{\underset{F,\varepsilon, n}{\rightharpoonup}} \{y' \} \}$ is called a {\bf ${\rm(}\bm{\varepsilon}$-{\rm)}terminal state} of $y \overset{\forall}{\underset{F,\varepsilon, n}{\rightharpoonup}} A$.
\end{definition}

It should be noted that the similar arrows $y \overset{\forall}{\underset{F,\varepsilon, n}{\rightharpoonup}} A$ and $y \overset{\forall}{\underset{F,\varepsilon}{\rightharpoondown}} A$ represent distinct symbols.
Notice that any dead-end terminal state of $y \overset{\forall}{\underset{F,\varepsilon}{\rightharpoondown}} A$ lies in the complement of the domain of $F$ and thus can no longer be mapped, whereas terminal states of $y \overset{\forall}{\underset{F,\varepsilon,n}{\rightharpoonup}} A$ are allowed to remain within the domain.

\begin{definition}
Let $\varepsilon \in (0,\infty]$ be a number. 
For any point $y \in Y$ and any subset $A \subseteq Y$ and for any $n \in \Z_{\geq 0}$, the relation $y \overset{\forall}{\underset{F,\varepsilon_{\Sigma}, n}{\rightharpoonup}} A$ holds if $\emptyset \neq \{ y' \mid y \overset{\exists}{\underset{F,\varepsilon_{\Sigma}, n}{\rightharpoonup}} \{y' \} \} \subseteq A$. 
Then any element of $\{ y' \mid y \overset{\exists}{\underset{F,\varepsilon_{\Sigma}, n}{\rightharpoonup}} \{y' \} \}$ is called a {\bf  ${\rm(}\bm{\varepsilon_{\Sigma}}$-{\rm)}terminal state} of $y \overset{\forall}{\underset{F,\varepsilon_{\Sigma}, n}{\rightharpoonup}} A$.
\end{definition}

Notice that any dead-end terminal state of $y \overset{\forall}{\underset{F,\varepsilon_{\Sigma}}{\rightharpoondown}} A$ lies in the complement of the domain of $F$ and thus can no longer be mapped, whereas terminal states of $y \overset{\forall}{\underset{F,\varepsilon_{\Sigma},n}{\rightharpoonup}} A$ are allowed to remain within the domain.

\section{Time extension of partial maps}
\label{sec:time_extension}

In this section, we show the following result. 
When the terminal state belongs to either $G$ or $B$, for any $\varepsilon \in [0,\infty)$, the whole space can be decomposed into two disjoint subsets: one consisting of elements that can be $\varepsilon$-controlled to reach $G$, and the other consisting of elements that only end up in $B$ even if $\varepsilon$-controlled.
This observation further implies that the set of desired states can be determined not only by direct computation, but also indirectly by identifying their complementary states.

Moreover, for any $\varepsilon \in \R$, there is a corresponding subset of elements that can be steered to certain dead-end terminal states under $\varepsilon$-control. 
These subsets effectively partition the whole space and thereby induce a filtration, which can be used in data analysis and is, in fact, applied in the next section.

\subsection{Setting of time extension of partial maps}

For any $a < b \in \Z$, put $[a,b]_{\Z} := [a,b] \cap \Z$. 
For a partial map $f \colon X \rightharpoonup X$ and any $a < b \in \Z_{\geq 0}$, put $Y := [a,b]_{\Z} \times X$ and we define a partial map $\widetilde{f} \colon Y \rightharpoonup Y$ by $\widetilde{f}(t,x) := (t+1, f(x))$. 
Then $\widetilde{f}$ is called the {\bf time extension} of $f$. 
Notice that $\mathop{\mathrm{dom}}\widetilde{f} \subseteq [a,b-1]_{\Z} \times X$ and $\mathop{\mathrm{Im}}\widetilde{f} \subseteq [a+1,b]_{\Z} \times X$.
For the time extension $\widetilde{f} \colon Y \rightharpoonup Y$, define a subset $D_{F} \subseteq Y$ as follows: 
\[
D_{F} := \{ (b-n,x) \in Y \mid n \in \{ 0,1, \ldots, b-a\}, f^n(x) \in X \}
\]
Put $F := \widetilde{f} \vert_{D_{F}} \colon D_{F} \rightharpoonup D_{F}$. 
Then $D_F - \mathop{\mathrm{dom}}F = \{b \} \times X = D_F \cap (\{b \} \times X)$. 
Notice that we introduce $D_F$ to justify the setting in which we consider what happens at the final time $b$ using $\mathop{\mathrm{dom}} F$, and also because a suitable filtration exists on $D_F$.

Suppose that $X$ has a metric $d$. 
Then define a non-degenerate cost function $c_d \colon D_{F} \times D_{F} \to [0,\infty]$, called the {\bf cost function associated to $\bm{d}$}, as follows: 
\[
c_d ((t,x),(t',x')) := \begin{cases}
d(x,x') & \text{if } t=t'\\
\infty & \text{if } t \neq t' 
\end{cases}
\]

From now on, in this section, suppose that $X$ has a metric $d$ and that $D_F$ equips the non-degenerate cost function $c_d$ associated with $d$. 
We have the following equivalences. 

\begin{lemma}
Let $A = \{b\} \times A_X \subseteq \{ b\} \times X = D_{F} - \mathop{\mathrm{dom}}F$ be a subset.
Then we have the following equalities for any $\varepsilon \in [0, \infty)$:
\[
\begin{split}
A_{F,\varepsilon} & = \{ (b-n,x) \in D_{F} \mid (b-n,x) \overset{\exists}{\underset{F,\varepsilon,n}{\rightharpoonup}} A \}
\\
& = \{ (b-n,x) \in D_{F} \mid x \overset{\exists}{\underset{f,\varepsilon,n}{\rightharpoonup}} A_X \}   
\end{split}
\]
\[
\begin{split}
A_{F,0} = A_{F,-0} &= \{ (b-n,x) \in D_{F} \mid f^n(x) \in A_X \}
\\
&= \{ y \in D_F \mid O^{\geq 0}_F(y) \cap A \neq \emptyset \}
\end{split}
\]
\[
\begin{split}
A_{F,-\varepsilon} &= \{ (b-n,x) \in D_{F} \mid (b-n,x) \overset{\forall}{\underset{F,\varepsilon,n}{\rightharpoonup}} A \} 
\\
&= \{ (b-n,x) \in D_{F} \mid x \overset{\forall}{\underset{f,\varepsilon,n}{\rightharpoonup}} A_X \} 
\end{split}
\]
\end{lemma}

\begin{proof}
By definition of $A_{F,\pm\varepsilon}$, we have the following equality for any $\varepsilon \in [0, \infty)$: 
\[
\begin{split}
A_{F,\varepsilon} & = \{ y \in D_{F} \mid y \overset{\exists}{\underset{F,\varepsilon}{\rightharpoonup}} A \} 
\\
& = \{ (b-n,x) \in D_{F} \mid (b-n,x) \overset{\exists}{\underset{F,\varepsilon}{\rightharpoonup}} A \}
\\
& = \{ (b-n,x) \in D_{F} \mid (b-n,x) \overset{\exists}{\underset{F,\varepsilon,n}{\rightharpoonup}} A \}
\\
&= \{ (b-n,x) \in D_{F} \mid x \overset{\exists}{\underset{f,\varepsilon,n}{\rightharpoonup}} A_X \} 
\end{split}
\] 
Moreover, we have the following equivalence: 
\[
\begin{split}
A_{F,0} & = \{ (b-n,x) \in D_{F} \mid x \overset{\exists}{\underset{f,0,n}{\rightharpoonup}} A_X \} 
= \{ (b-n,x) \in D_{F} \mid f^n(x) \in A_X \} 
\end{split}
\] 
For any $y = (b-n,x) \in A_{F,0}$, we have $[y]_F^{0} \setminus \mathop{\mathrm{dom}}F = \{ (b,f^n(x)) \} \subseteq A \subseteq D_{F} - \mathop{\mathrm{dom}}F$. 
Then we obtain $A_{F,0} = \{ y \in D_F \mid O^{\geq 0}_F(y) \cap A \neq \emptyset \}$. 
Furthermore, we have the following equalities:
\[
\begin{split}
A_{F,-0} & = \{ y \in A_{F,0} \mid \emptyset \neq [y]_F^{0} \setminus \mathop{\mathrm{dom}}F \subseteq A \} 
\\
& = \{ (b-n,x) \in A_{F,0} \mid \{ (b,f^n(x)) \} \subseteq A \} 
\\
& = \{ (b-n,x) \in A_{F,0} \mid f^n(x) \in A_X \}
\\
& = \{ (b-n,x) \in D_{F} \mid f^n(x) \in A_X \} = A_{F,0}
\\
A_{F,-\varepsilon} &= \{ y \in A_{F,-0} \mid y \overset{\forall}{\underset{F,\varepsilon}{\rightharpoondown}} A \} 
\\
&= \{ y \in A_{F,-0} \mid \emptyset \neq [y]_F^\varepsilon \setminus \mathop{\mathrm{dom}}F \subseteq A \} 
\\
&= \{ (b-n,x) \in A_{F,-0} \mid \emptyset \neq [(b-n,x)]_F^\varepsilon \setminus \mathop{\mathrm{dom}}F \subseteq A \} 
\\
&= \{ (b-n,x) \in D_{F} \mid \emptyset \neq [(b-n,x)]_F^\varepsilon \setminus \mathop{\mathrm{dom}}F \subseteq A, 
f^n(x) \in A_X \} 
\end{split}
\]

By $D_F - \mathop{\mathrm{dom}}F = \{b \} \times X$ and $D_{F} = \{ (b-n,x) \in Y \mid n \in \{ 0,1, \ldots, b-a\}, f^n(x) \in X \}$, 
the following equivalence holds for any $(b-n,x) \in D_{F}$:
\[
\begin{split}
& \emptyset \neq [(b-n,x)]_F^\varepsilon \setminus \mathop{\mathrm{dom}}F \subseteq A, f^n(x) \in A_X
\\
\Longleftrightarrow \,\, & \emptyset \neq [(b-n,x)]_F^\varepsilon \setminus \mathop{\mathrm{dom}}F \subseteq A
\\
\Longleftrightarrow \,\, & \emptyset \neq \{ y' \in D_F  \mid (b-n,x) \overset{\exists}{\underset{F,\varepsilon}{\rightharpoonup}} \{y' \} \} \setminus \mathop{\mathrm{dom}}F \subseteq A
\\
\Longleftrightarrow \,\, & \emptyset \neq \{ y' \in D_F - \mathop{\mathrm{dom}}F \mid (b-n,x) \overset{\exists}{\underset{F,\varepsilon}{\rightharpoonup}} \{y' \} \} \subseteq A
\\
\Longleftrightarrow \,\, & \emptyset \neq \{b \} \times \{ x' \in X \mid (b-n,x) \overset{\exists}{\underset{F,\varepsilon}{\rightharpoonup}} \{(b,x') \} \} \subseteq  \{b\} \times A_X = A
\\
\Longleftrightarrow \,\, & \emptyset \neq \{ x' \in X \mid (b-n,x) \overset{\exists}{\underset{F,\varepsilon,n}{\rightharpoonup}} \{(b,x') \} \} \subseteq A_X
\\
\Longleftrightarrow \,\, & \emptyset \neq \{ x' \in X \mid x \overset{\exists}{\underset{f,\varepsilon,n}{\rightharpoonup}} \{x' \} \} \subseteq A_X
\\
\Longleftrightarrow \,\, & x \overset{\forall}{\underset{f,\varepsilon,n}{\rightharpoonup}} A_X
\\
\Longleftrightarrow \,\, & (b-n,x) \overset{\forall}{\underset{F,\varepsilon,n}{\rightharpoonup}} A
\end{split}
\]
This implies $A_{F,-\varepsilon} = \{ (b-n,x) \in D_{F} \mid x \overset{\forall}{\underset{f,\varepsilon,n}{\rightharpoonup}} A_X \} = \{ (b-n,x) \in D_{F} \mid (b-n,x) \overset{\forall}{\underset{F,\varepsilon,n}{\rightharpoonup}} A \}$. 
\end{proof}

By the same argument of the previous proof, we have the following statement. 

\begin{lemma}
Let $A = \{b\} \times A_X \subseteq \{ b\} \times X = D_{F} - \mathop{\mathrm{dom}}F$ be a subset.
Then we have the following equalities for any $\varepsilon \in [0, \infty)$:
\[
\begin{split}
A_{F,\varepsilon_{\Sigma}} & = \{ (b-n,x) \in D_{F} \mid (b-n,x) \overset{\exists}{\underset{F,\varepsilon_{\Sigma},n}{\rightharpoonup}} A \} 
\\
& = \{ (b-n,x) \in D_{F} \mid x \overset{\exists}{\underset{f,\varepsilon_{\Sigma},n}{\rightharpoonup}} A_X \} 
\end{split}
\]
\[
\begin{split}
A_{F,0_{\Sigma}} & = A_{F,-0_{\Sigma}} = \{ (b-n,x) \in D_{F} \mid f^n(x) \in A_X \} 
\\
& = \{ y \in D_F \mid O^{\geq 0}_F(y) \cap A \neq \emptyset \} 
\end{split}
\]
\[
\begin{split}
A_{F,-\varepsilon_{\Sigma}} & = \{ (b-n,x) \in D_{F} \mid (b-n,x) \overset{\forall}{\underset{F,\varepsilon_{\Sigma},n}{\rightharpoonup}} A \} 
\\
& = \{ (b-n,x) \in D_{F} \mid x \overset{\forall}{\underset{f,\varepsilon_{\Sigma},n}{\rightharpoonup}} A_X \} 
\end{split}
\]
\end{lemma}

\subsection{Dichotomy for final states}
We show the following dichotomy for final states.

\begin{theorem}\label{thm:eps_GBseparation}
Let $G,B \subset \{b \} \times X = D_{F} - \mathop{\mathrm{dom}}F$ be subsets with $G \sqcup B = \{b \} \times X$. 
Suppose that $X$ has a metric $d$ and that $D_F$ equips the cost function associated with $d$.
For any $\varepsilon \in [0,\infty)$, we have $G_{F,\varepsilon} \sqcup B_{F,-\varepsilon} = D_F$. 
\end{theorem}

\begin{proof}
Since $D_{F} - \mathop{\mathrm{dom}}F = \{b \} \times X$ and $D_{F} = \{ (b-n,x) \in Y \mid n \in \{ 0,1, \ldots, b-a\}, f^n(x) \in X \}$, we have the following equalities for any $(b-n,x) \in D_{F}$: 
\[
\begin{split}
&~(b-n,x) \in B_{F,-\varepsilon} 
\\
\Longleftrightarrow&~(b-n,x) \overset{\forall}{\underset{F,\varepsilon,n}{\rightharpoonup}} B
\\
\Longleftrightarrow&~\emptyset \neq \{ y' \in D_{F} - \mathop{\mathrm{dom}}F \mid (b-n,x) \overset{\exists}{\underset{F,\varepsilon,n}{\rightharpoonup}} \{y' \} \} \subseteq B
\\
\Longleftrightarrow&~\{ y' \in D_{F} - \mathop{\mathrm{dom}}F \mid (b-n,x) \overset{\exists}{\underset{F,\varepsilon,n}{\rightharpoonup}} \{y' \} \} \subseteq B
\\
\Longleftrightarrow&~\{ y' \in D_{F} - \mathop{\mathrm{dom}}F \mid (b-n,x) \overset{\exists}{\underset{F,\varepsilon}{\rightharpoonup}} \{y' \} \} 
\subseteq B = (\{b \} \times X) - G
\\
\Longleftrightarrow&~\{ (b,x') \in \{b \} \times X \mid (b-n,x) \overset{\exists}{\underset{F,\varepsilon}{\rightharpoonup}} \{(b,x') \} \} 
\subseteq B = (\{b \} \times X) - G
\\
\Longleftrightarrow&~\emptyset = \{ (b,x') \in \{b \} \times X \mid (b-n,x) \overset{\exists}{\underset{F,\varepsilon}{\rightharpoonup}} \{(b,x') \} \} \cap G
\\
\Longleftrightarrow&~\emptyset = \{ y' \in D_F \mid (b-n,x) \overset{\exists}{\underset{F,\varepsilon}{\rightharpoonup}} \{y' \} \} \cap G
\\
\Longleftrightarrow&~\emptyset = [(b-n,x)]^\varepsilon_F \cap G
\\
\Longleftrightarrow&~(b-n,x) \notin G_{F,\varepsilon} 
\end{split}
\]
This implies the assertion.
\end{proof}

By the same argument of the previous proof, we have the following statement. 

\begin{theorem}\label{thm:eps_sum_GBseparation}
Let $G,B \subset \{b \} \times X = D_{F} - \mathop{\mathrm{dom}}F$ be subsets with $G \sqcup B = \{b \} \times X$. 
Suppose that $X$ has a metric $d$ and that $D_F$ equips the cost function associated with $d$.
For any $\varepsilon \in [0,\infty)$, we have $G_{F,\varepsilon_{\Sigma}} \sqcup B_{F,-\varepsilon_{\Sigma}} = D_F$. 
\end{theorem}

The conditions $G,B \subset D_{F} - \mathop{\mathrm{dom}}F$ in the previous two theorems are necessary. 
In fact, there is a map $f \colon \{0,1,2\} \to \{0,1,2\}$ defined by $0 = f(0)$ and $1 = f(1) = f(2)$, and define a cost function $c \colon \{0,1,2\}^2 \to \{0,1,2\}$ by $c(n,m) := \vert n - m \vert$. 
Set $G := \{ 2 \}$ and $B := \{0 \}$. 
Then $B_{F,-1/2} = B_{F,-1/2_{\sum}} = \{ 0 \}$ and $G_{F,1/2} = G_{F,1/2_{\sum}} = \{ 2 \} \subsetneq \{1,2 \} = X - B_{F,-1/2}$. 
Moreover, $B_{F,-(1+\varepsilon)} = \emptyset$ and $G_{F,\varepsilon} = \{2 \}$ for any $\varepsilon \in [0,\infty)$. 


\begin{figure*}[ht]
  \centering
  \includegraphics[width=\textwidth]{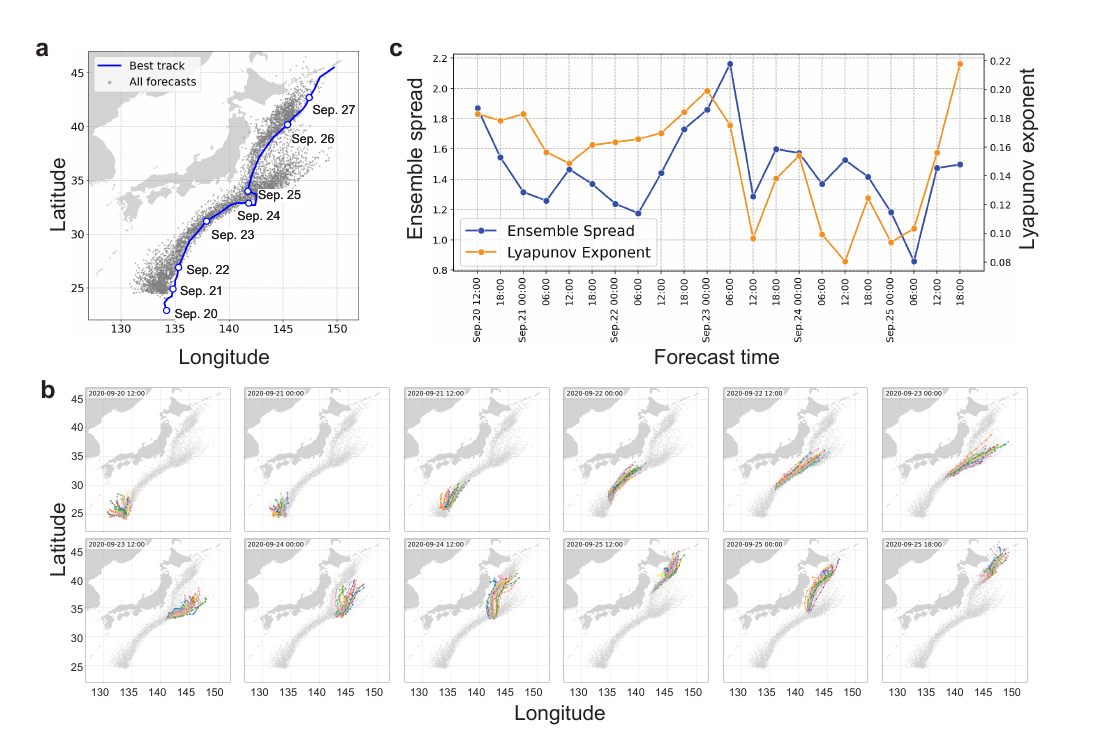}
  \caption{
    Ensemble forecast dataset for tropical cyclone Dolphin.
    \textbf{a}, Best‐track trajectory (real orbit, blue line) overlaid with all forecast center positions (gray dots) of Dolphin.
    \textbf{b}, Selected ensemble forecast tracks at representative initial times: each panel shows 21 ensemble members (colored lines), plotted against the background of all forecasts (gray dots). 
    The times correspond to the initial times of the forecasts.
    {
    \textbf{c}, Time evolution of ensemble spread (blue) and a finite-time Lyapunov exponent proxy (orange). The ensemble spread is computed as the spatial variance of ensemble members at each forecast time, while the Lyapunov exponent is estimated from the growth rate of pairwise distances between trajectories. The comparison illustrates how perturbation growth and forecast dispersion evolve over time.
    }
    }
  \label{fig:application_dataset}
\end{figure*}

\subsection{Filtrations on $D_F$}

The following statement ensures that, for two disjoint terminal states, a filtration exists on the complement of one terminal state with respect to the other, which is applied in the next section.

\begin{theorem}\label{thm:eps_GBcovering}
Let $G,B \subset \{b \} \times X = D_{F} - \mathop{\mathrm{dom}}F$ be subsets with $G \sqcup B = \{b \} \times X$. 
Suppose that $X$ has a metric $d$ and that $D_F$ equips the cost function associated with $d$.
The following statements hold: 
\\
{\rm(1)} If there is a point $(b-n,x) \in \mathop{\mathrm{dom}}F$ with $F^n(b-n,x)\in G$ for some positive integer $n$, then $D_F - B = \bigcup_{\varepsilon \in \R} G_{F,\varepsilon_{\Sigma}} = \bigcup_{\varepsilon \in \R} G_{F,\varepsilon}$. 
\\
{\rm(2)} If there is a point $(b-n,x) \in \mathop{\mathrm{dom}}F$ with $F^n(b-n,x)\in B$ for some positive integer $n$, then $D_F - G = \bigcup_{\varepsilon \in \R} B_{F,\varepsilon_{\Sigma}} = \bigcup_{\varepsilon \in \R} B_{F,\varepsilon}$. 
\end{theorem}

\begin{proof}
By symmetry, it suffices to show assertion (1). 
Fix a point $p = (b-n,x) \in \mathop{\mathrm{dom}}F$ with $F^n(b-n,x)\in G$ for some positive integer $n$. 
Choose any point $q = (b-n',x') \in \mathop{\mathrm{dom}}F$ for some positive integer $n'$. 
Put $\varepsilon' := d(f^{n'-1}(x'),f^{n-1}(x)) \in [0,\infty)$.  
From $F(F^{n-1}(b-n,x)) = F^n(b-n,x)\in G$, the sequence 
\[
(q_0,q_1, \ldots, q_{n'-2},q_{n'-1}) := (q,F(q), \ldots , F^{n'-2}(q), F^{n-1}(p))
\]
satisfies  $c_d(q,q_0) = 0$, $c_d(F(q_i),q_{i+1})= 0$ for any $i \in \{1, \ldots , n'-3\}$,  
\[
c_d(F(q_{n'-2}),q_{n'-1}) = c_d(F^{n'-1}(q),F^{n-1}(p)) = \varepsilon'
\]
 and $F(q_{n'-1}) = F^n(b-n,x)\in G$. 
This means that $q \overset{\exists}{\underset{F,\varepsilon'_{\Sigma}}{\rightharpoonup}} G$ and so $q \overset{\exists}{\underset{F,\varepsilon'}{\rightharpoonup}} G$. 
Therefore, we obtain $\bigcup_{\varepsilon \in \R} G_{F,\varepsilon_{\Sigma}} \supseteq G \sqcup \mathop{\mathrm{dom}}F = D_F - B$. 
By $G_{F,\varepsilon_{\Sigma}} \subseteq G_{F,\varepsilon} \subseteq G \sqcup \mathop{\mathrm{dom}}F = D_F - B$, assertion (1) holds. 
\end{proof}

\section{Application to Meteorology}

This section applies the proposed framework to ensemble weather forecast data for tropical cyclone (typhoon) ``Dolphin'', generated by the Meso-scale Ensemble Prediction System (MEPS)\cite{oettli2024objective,oettli2025meteorological}. 
Dolphin formed in September 2020 and was initially forecast to track eastward over the western Pacific to the east of Japan; however, its trajectory shifted to a northward route, causing it to graze the Tohoku region and deliver heavy rainfall across much of Japan (Figure~\ref{fig:application_dataset}a--b). 
In other words, the forecast data for Dolphin exhibit two distinct pathways: eastward and northward, which correspond to a “good” cluster and a “bad” cluster, respectively. 
The objective of this section is to quantify how perturbations of magnitude at most $\varepsilon$ affect the eventual outcome of the trajectory.
The dataset used here was provided by the Japan Meteorological Agency (JMA), and the code for reproducing our analysis is available at \url{https://github.com/yusuke-imoto-lab/eps-attracting-basin}, and Zendo \cite{ImotoYokoyama2026_eps_attracting_basin}.

\begin{figure*}[h]
  \centering
  \includegraphics[width=\textwidth]{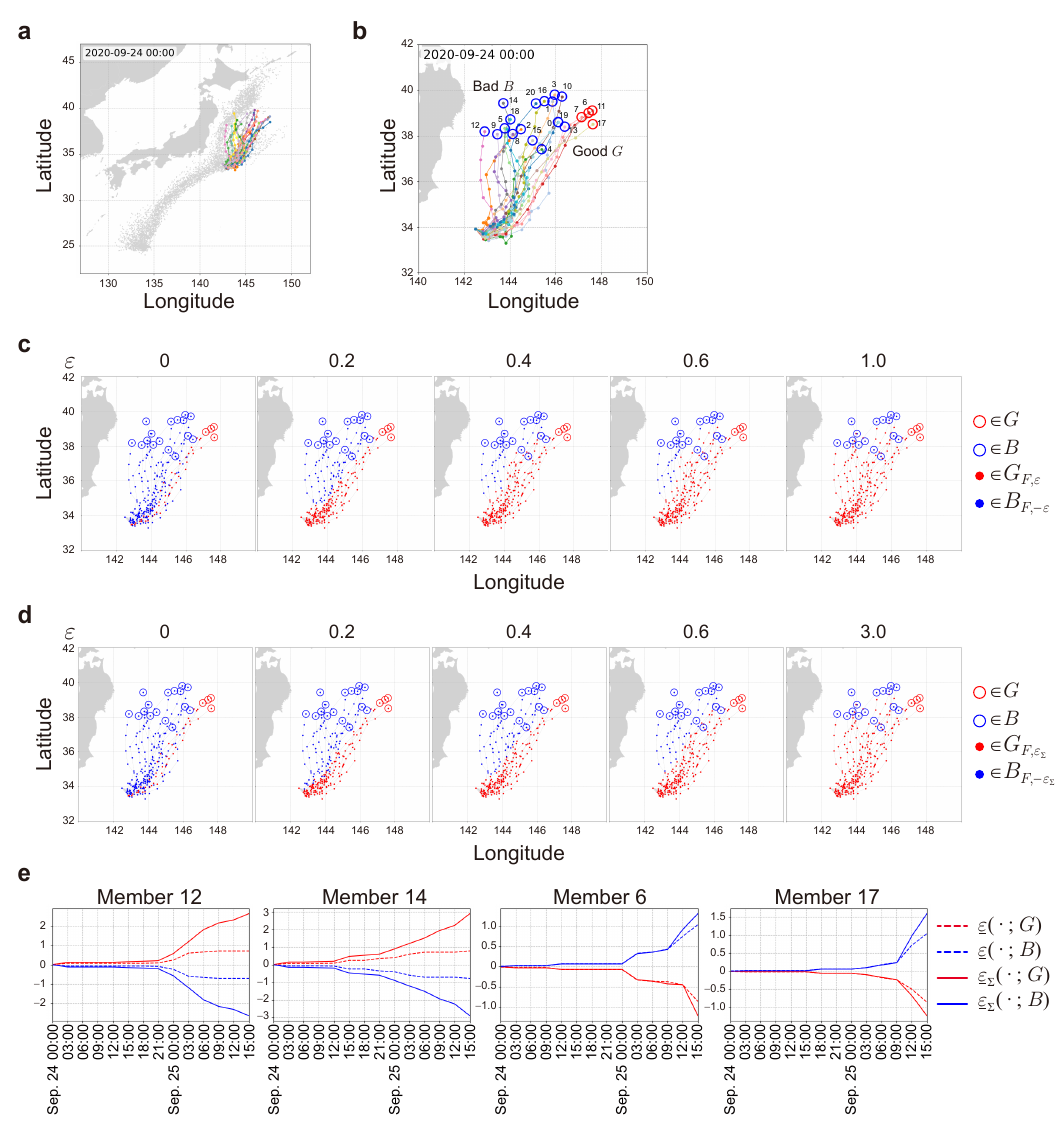}
  \caption{
    Numerical experiment using one ensemble weather forecast dataset.   
    \textbf{a}, Ensemble forecast tracks of the tropical cyclone initialized at 00:00 UTC on September 24, 2020, showing a wide east–west spread of possible trajectories. 
    \textbf{b}, Partition of final‐time ensemble points into the good cluster $G$ (eastward tracks, red circles) and the bad cluster $B$ (northward tracks, blue circles).  
    \textbf{c}--\textbf{d}, $\varepsilon$‐ and $\varepsilon_{\Sigma}$‐attracting basins for $G$ (red) and $B$ (blue) at various values of $\varepsilon$.
    \textbf{e}, Time series of $\underline{\varepsilon}(\,\cdot\,; G)$ and $\underline{\varepsilon}(\,\cdot\,; B)$ (solid red and blue) and of $\underline{\varepsilon}_{\Sigma}(\,\cdot\,; G)$ and $\underline{\varepsilon}_{\Sigma}(\,\cdot\,; B)$ (dashed red and blue) for four representative ensemble members (the member numbers corresponding to \textbf{b}). 
  } \label{fig:application_one_data}
\end{figure*}

The dataset consists of six-hourly forecasts of the Dolphin's center position from 12:00 UTC on September 20, 2020 to 18:00 UTC on September 25, 2020. 
At each six-hour interval, the initial state was newly set based on weather observations, and 21 ensemble forecasts were generated using the small-perturbed initial states based on a data assimilation system. 
Each ensemble forecast is a time series consisting of 13 timepoints with three-hour time intervals, covering a total forecast length of 39 hours (Figure~\ref{fig:application_dataset}a--b).

{
In addition, we compute standard indicators of forecast sensitivity for comparison. 
Figure~\ref{fig:application_dataset}c shows the time evolution of the ensemble spread and a finite-time Lyapunov exponent proxy estimated from the growth rate of pairwise distances between ensemble trajectories. 
The ensemble spread increases around 00:00 UTC on September 23 and remains relatively high until around 06:00 UTC, after which it decreases and stays comparatively low for a while before increasing again near the end of the period. 
The Lyapunov-based quantity shows a similar temporary increase before dropping to lower values around September 23 to September 24, followed by another rise toward the end of the forecast window. 
These standard indicators indeed capture changes in forecast sensitivity at certain times. 
However, they do not directly indicate whether the forecast uncertainty is associated with transitions between the good and bad clusters. 
This motivates the use of our framework, which is designed to quantify sensitivity with respect to such qualitatively distinct forecast outcomes.
}

\subsection{Validation using one ensemble weather forecast dataset. }
First, we validate the proposed method using the ensemble weather forecast dataset initialized at 00:00 UTC on September 24, 2020 (Figure~\ref{fig:application_one_data}a).  
At this initial time, the tropical cyclone forecast tracks were spread widely in the east-west direction, making it a critical time to determine the eventual path of the tropical cyclone \cite{oettli2024objective}.  
The dataset consists of 21 ensemble members over 14 forecast times (initial time and 13 forecast steps); we denote by $X\subset\mathbb{R}^2$ all its center positions ($|X| = 21 \times 14 = 294$), and set $Y =  [0,13]_{\mathbb{Z}}\times X$, which corresponds to the setting in Section \ref{sec:time_extension}. 
We set the metric $d$ on $X\subset\mathbb{R}^2$ to be the Euclidean distance.

Among these 21 ensemble members, we classified the four tracks that move eastward at the final forecast time as the good cluster $G$, and the remaining tracks as the bad cluster $B$ (Figure~\ref{fig:application_one_data}b). 
For each of these clusters, we computed the attracting basins $G_{F,\varepsilon}, B_{F,-\varepsilon}, G_{F,\varepsilon_{\Sigma}}, B_{F,-\varepsilon_{\Sigma}}$ (Figure~\ref{fig:application_one_data}c--d). 
Corresponding to Theorem \ref{thm:eps_GBseparation}, for every $\varepsilon\ge0$, each point in $Y$ belongs to exactly one of $G_{F,\varepsilon}$ or $B_{F,-\varepsilon}$.  
Moreover, reflecting Theorem \ref{thm:eps_GBcovering}, as $\varepsilon$ increases, all points except those in $B$ eventually lie in $G_{F,\varepsilon}$. 
The same conclusions hold for $G_{F,\varepsilon_{\Sigma}}$ and $B_{F,-\varepsilon_{\Sigma}}$, following Theorems \ref{thm:eps_sum_GBseparation} and \ref{thm:eps_GBcovering}.  

Comparing the $\varepsilon$‐ and $\varepsilon_{\Sigma}$‐attracting basins, we observe that the two distributions coincide at $\varepsilon = 0$; and the inclusions $G_{F,\varepsilon_{\Sigma}} \;\subseteq\; G_{F,\varepsilon}$ and $B_{F,-\varepsilon} \;\subseteq\; B_{F,-\varepsilon_{\Sigma}}$
hold for $\varepsilon>0$, consistent with Lemmas~\ref{lem:two_inclusions} and \ref{lem:inclusions}.
{
In other words, the region where perturbations of magnitude at most $\varepsilon$ at each time step can influence the eventual outcome (Figure~\ref{fig:application_one_data}c) is larger than the corresponding region when the total perturbation magnitude over all times is bounded by $\varepsilon$ (Figure~\ref{fig:application_one_data}d), which agrees with our intuition.
}

In addition, we computed the time evolution of the debut functions $\underline{\varepsilon}$ and $\underline{\varepsilon}_{\Sigma}$ with respect to clusters $G$ and $B$ for each ensemble member (Figure~\ref{fig:application_one_data}e). 
At every forecast time, the signs of $\underline{\varepsilon}$ and $\underline{\varepsilon}_{\Sigma}$ for the same cluster coincide; and the absolute values of the debut functions $\underline{\varepsilon}$ are always equivalent to or less than those of $\underline{\varepsilon}_{\Sigma}$, in accordance with Lemma~\ref{lem:debut_ineq1}. 
Furthermore, the debut functions $\underline{\varepsilon}$ and $\underline{\varepsilon}_{\Sigma}$ are either non-increasing or non-decreasing in time, reflecting Lemma~\ref{lem:debut_ineq2}.  

\subsection{Detection of controllable region via debut functions.}
\begin{figure*}[pht]
  \centering
  \includegraphics[width=\textwidth]{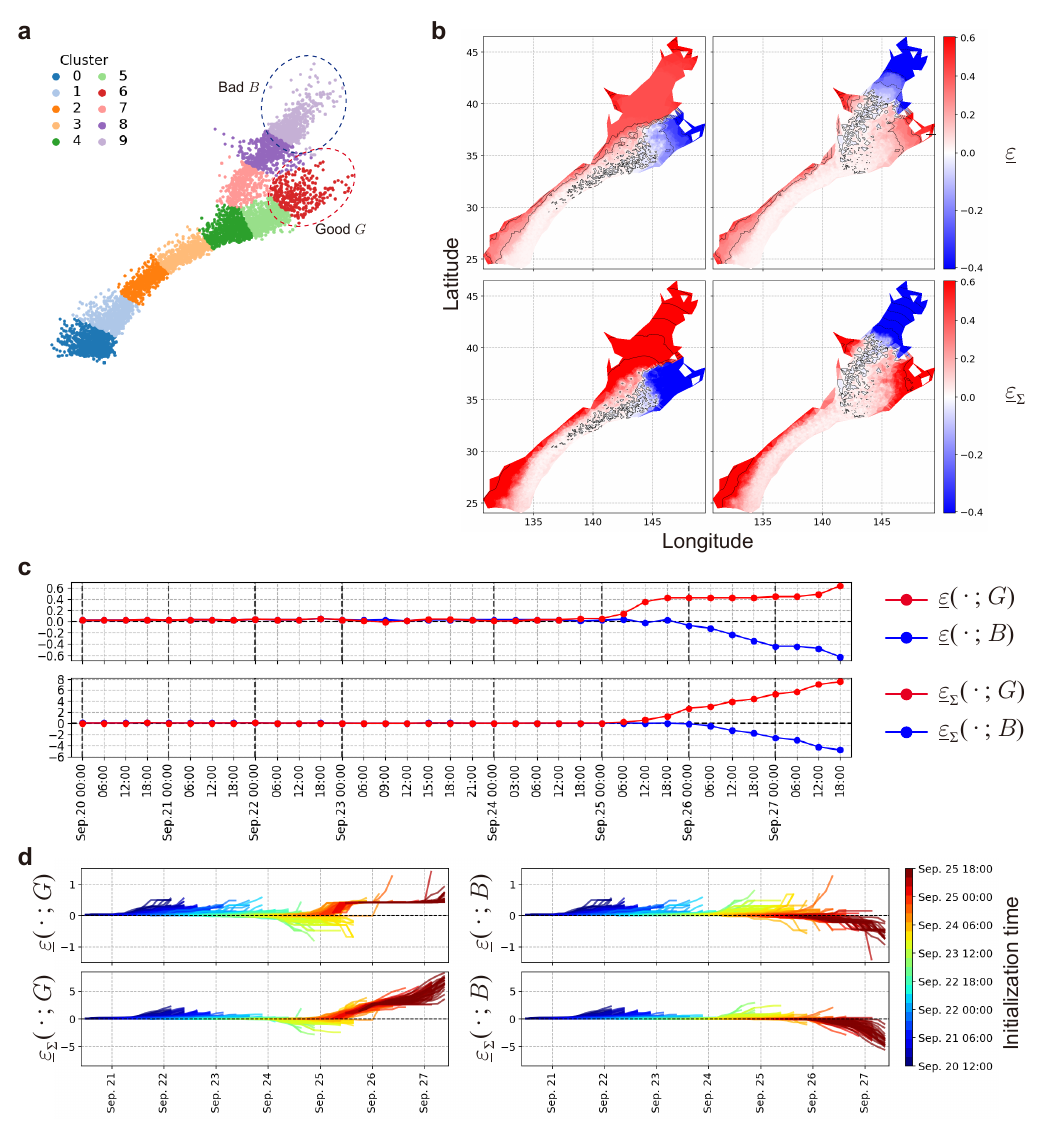}
  \caption{
    Numerical experiment using the full ensemble weather forecast dataset. 
    \textbf{a}, $k$-means clustering of all the forecast center positions into 10 clusters. 
    Clusters 6 and 9 were assigned as good cluster $G$ and bad cluster $B$, respectively. 
    \textbf{b}, Distribution of $\underline{\varepsilon}$ (top) and $\underline{\varepsilon}_{\Sigma}$ (bottom) for clusters $G$ (reft) and $B$ (right). 
    \textbf{c}, Time evolution of $\underline{\varepsilon}$ (top) and $\underline{\varepsilon}_{\Sigma}$ (bottom) for the good cluster $G$ (red) and the bad cluster $B$ (blue) on the best track.  
    {
    \textbf{d}, Time evolution of the debut functions for each ensemble member, shown separately for the good cluster $G$ (left) and the bad cluster $B$ (right). The top panels show $\underline{\varepsilon}$ and the bottom panels show $\underline{\varepsilon}_{\Sigma}$. The color indicates the initialization time of each ensemble forecast. 
    }
    }
  \label{fig:application_all_data}
\end{figure*}

Next, we conduct a numerical experiment using the full ensemble weather forecast dataset of Dolphin’s center positions. 
Since each ensemble forecast was computed independently at each initial time, the time stamps of individual tracks are not synchronized. 
Therefore, we set $Y\subset\mathbb{R}^2$ the two‐dimensional point cloud of Dolphin’s center positions, which corresponds to the setting in Section \ref{sec:Preliminaries}. 
We set the cost function $c$ to be the Euclidean distance. 
We applied K-means clustering to the point cloud $Y$ and assigned the east and north clusters at the final status as the good cluster $G$ and bad cluster $B$, respectively (Figure~\ref{fig:application_all_data}a).
Next, we computed the debut functions $\underline{\varepsilon}$ and $\underline{\varepsilon}_{\Sigma}$ with respect to a subset $G$ and $B$ and visualized them (Figure~\ref{fig:application_all_data}b).
Because the $\varepsilon$- and $\varepsilon_{\Sigma}$-attracting basins are filtrations (Theorem~\ref{thm:filtration}), $A_{F,\varepsilon}$ and $A_{F,\varepsilon_\Sigma}$ coincide with the sublevel set of $\underline{\varepsilon}(\,\cdot\,; A)$ and $\underline{\varepsilon}_{\Sigma}(\,\cdot\,; A)$, respectively. 

The spatial distributions of $\underline{\varepsilon}$ and $\underline{\varepsilon}_{\Sigma}$ exhibit several notable features. 
West of longitude $140^\circ$E, both $\underline{\varepsilon}$ and $\underline{\varepsilon}_{\Sigma}$ are positive almost everywhere, indicating that trajectories in these regions are not automatically drawn into either the good cluster $G$ or the bad cluster $B$ without control. 
Consequently, these areas represent high-uncertainty zones in which small perturbations could steer the cyclone toward either cluster under an appropriate control input. 
Around the intersection of $140^\circ$E and $35^\circ$N, negative-valued regions of $\underline{\varepsilon}(\,\cdot\,; G)$ and $\underline{\varepsilon}(\,\cdot\,; B)$ (and similarly for $\underline{\varepsilon}_{\Sigma}$) are intermingled, bringing to mind a riddled basin or a fractal basin boundary \cite{alexander1992riddled, mcdonald1985fractal}. 
Further, northeast, near $145^\circ$E and $35^\circ$N, $\underline{\varepsilon}(\,\cdot\,; G)$ becomes negative while $\underline{\varepsilon}(\,\cdot\,; B)$ remains positive, indicating that points in this region will naturally converge to the good cluster. 
{
This implies that, for the trajectory to evolve toward the good cluster, perturbations must shift the center position into regions where $\underline{\varepsilon}(\,\cdot\,;G)<0$ and $\underline{\varepsilon}(\,\cdot\,;B)>0$ (or where $\underline{\varepsilon}_\Sigma(\,\cdot\,;G)<0$ and $\underline{\varepsilon}_\Sigma(\,\cdot\,;B)>0$).
}

In addition, we computed $\underline{\varepsilon}$ and \(\underline{\varepsilon}_{\Sigma}\) along the best track by referring to the values on the nearest forecast point from each best track point (Figure~\ref{fig:application_all_data}c). 
Until 00:00 UTC on September 25, both $\underline{\varepsilon}(\,\cdot\,; G)$ and $\underline{\varepsilon}(\,\cdot\,; B)$ keep just above zero, indicating an indeterminate state in which a small perturbation could steer the cyclone toward either the good or bad cluster. 
After 00:00 UTC on September 25, however, $\underline{\varepsilon}(\,\cdot\,; G)$ increases while $\underline{\varepsilon}(\,\cdot\,; B)$ decreases; by 00:00 UTC on September 26, \(\underline{\varepsilon}(\,\cdot\,; B)<0\), implying that without control the trajectory of tropical cyclone is drawn inevitably into the bad cluster $B$. 
If controllable, for example, if we can apply a perturbation of magnitude $\varepsilon=0.2$ at each time step, then the tropical cyclone could still be diverted outside the bad cluster until 12:00 UTC on September 26. 
In that scenario, as we mentioned in the last paragraph, steering the system into $\underline{\varepsilon}(\,\cdot\,;G)$-negative and $\underline{\varepsilon}(\,\cdot\,;B)$-positive valued regions guarantees convergence to the good cluster.
When the total control energy is constrained (the \(\varepsilon_{\Sigma}\) case), the increase in $\underline{\varepsilon}_{\Sigma}(\,\cdot\,; G)$ and the decrease in $\underline{\varepsilon}_{\Sigma}(\,\cdot\,; B)$ are more pronounced, and the term during which effective control can be exerted under the same $\varepsilon$ is correspondingly shorter.

{
We finally compare our framework with standard indicators of forecast sensitivity. 
Figure~\ref{fig:application_dataset}c shows the time evolution of ensemble spread and a finite-time Lyapunov exponent proxy, while Figure~\ref{fig:application_all_data}d presents the debut functions computed for each ensemble member. 
The ensemble spread and Lyapunov-based measures capture the overall growth of perturbations and exhibit temporal variations, such as an increase around September 23 followed by a period of relatively low values. 
From these standard indicators alone, one might interpret that the eventual outcome is already determined around September 23. 
However, the forecast trajectories in Figure~\ref{fig:application_dataset}b show that the ensemble does not actually begin to move clearly toward the bad outcome until around September 24. 
This demonstrates that our framework complements standard measures by capturing outcome-level sensitivity, providing an additional perspective that is not directly accessible from them alone.

In contrast, the debut functions directly quantify the minimal perturbation required to induce transitions between the good and bad clusters. 
As shown in Figure~\ref{fig:application_all_data}d, ensemble members satisfying $\underline{\varepsilon}(\,\cdot\,;G)<0$ already appear around September 23, indicating the emergence of trajectories that are attracted to the good cluster. 
Subsequently, ensemble members satisfying $\underline{\varepsilon}(\,\cdot\,;B)<0$ begin to appear around September 24, reflecting the onset of trajectories that are attracted to the bad cluster. 
Thus, our framework captures the temporal structure of the transition more faithfully: it detects the early appearance of attraction toward the good cluster and the later emergence of attraction toward the bad cluster, rather than suggesting that the fate of the forecast is fixed at a single earlier time. 
This demonstrates that our framework captures outcome-level sensitivity, which is not directly accessible from standard measures.
}

\subsection{Discussion.}
The proposed framework is data-driven, with a range of theoretically guaranteed properties, as well as high computational flexibility and interpretability.
Because the input is multiple time series data generated by an underlying dynamical system $f$, computations can proceed even when the dynamical system $f$ itself is unknown. 
Moreover, since all calculations aside from the pairwise cost evaluations are independent of data dimensionality, even high-dimensional datasets do not lead to an extreme increase in computational cost. 
By leveraging the filtration property and visualizing sublevel sets (Figure~\ref{fig:application_all_data}b), one can intuitively and quantitatively identify the tipping points, namely “when” and “where” small perturbations can or cannot alter the eventual outcome. 
Finally, offering two formulations, $\varepsilon$ and $\varepsilon_{\Sigma}$, allows users to choose the most appropriate version depending on the nature of the phenomenon and the available perturbation models.

Because of these features, the proposed framework can be applied to high-dimensional datasets, such as large-scale meteorological data that include temperature, humidity, pressure, and wind speed, and it holds promise for applications to forecast uncertainty and sensitivity analysis in extreme weather events, including heavy rainfall. 
In particular, even when perturbations are restricted to certain variables, times, or locations, such constraints can be incorporated by appropriately redesigning the cost function. 

An important direction for future work is to extend the framework so as to account for reliability and uncertainty in a probabilistic setting. 
In practical applications, the estimation of sensitivity depends not only on the mathematical formulation but also on the quality of observational data, the accuracy of data assimilation, and the presence of process and observational noise. 
Therefore, combining our approach with uncertainty quantification methods and data assimilation techniques would enable a more robust assessment of forecast sensitivity in realistic high-dimensional systems. 
Going forward, we plan to develop such an integrated framework and to broaden the applicability of our method to these high-dimensional data.

Furthermore, while our current focus is on data-driven frameworks for systems with unknown dynamics, exploring a deterministic setting where the governing equations are explicitly available serves as an essential theoretical baseline. Such deterministic systems function as rigorous toy models to validate the mathematical consistency and reliability of the filtration framework before addressing more complex, real-world stochastic processes. A detailed analysis of these models, alongside their systematic extension into probabilistic contexts to account for observational noise and uncertainty, is also left for future work.

\subsection*{Acknowledgments}
This research was partially supported by the JST Moonshot R\&D Program and CREST (Grant Numbers JPMJMS2389 and JPMJCR24Q1). 
T.Y. acknowledges the hospitality of the Fields Institute for Research in Mathematical Sciences, where part of this work was revised during a stay.
The authors thank Dr. Takahito Mitsui and Dr. Keita Tokuda for their useful comments, and Dr. Shunji Kotsuki and Dr. Oettli Pascal for sharing ensemble weather forecast data.


\bibliographystyle{amsplain}
\bibliography{GBF}

\end{document}